\documentclass[10pt]{amsart}

\usepackage{verbatim}

\usepackage{amscd,amsmath,amssymb, hyperref}
\usepackage[mathscr]{eucal}
\usepackage{graphicx}

\theoremstyle{plain}

\numberwithin{equation}{section}

\newtheorem{theorem}{Theorem}[section]
\newtheorem{lemma}[theorem]{Lemma}
\newtheorem{proposition}[theorem]{Proposition}
\newtheorem{corollary}[theorem]{Corollary}
\newtheorem{question}[theorem]{Question}
\newtheorem*{question*}{Question}

\theoremstyle{definition}
\newtheorem{definition}{Definition}[section]

\theoremstyle{remark}
\newtheorem{remark}{Remark}[section]
\newtheorem*{claim*}{Claim}
\newtheorem*{example*}{Example}
\newtheorem{example}[remark]{Example}
\newtheorem*{remark*}{Remark}

\newcommand{\R}{\mathbf{R}}

\begin{document}

\title[cost-subdifferential]
{On the cost-subdifferentials of cost-convex functions
}

\author{Young-Heon Kim
and Robert J. McCann
}

\address{Department of Mathematics, University of Toronto\\
  Toronto, Ontario Canada M5S 2E4}

\email{yhkim@math.toronto.edu, mccann@math.toronto.edu
}
\date{\today}



\maketitle

\begin{abstract}
We are interested in the cost-convex potentials in optimal mass transport theory, and we show
by direct and geometric arguments the equivalence between cost-subdifferentials and
ordinary subdifferentials of cost-convex functions, under the assumptions {\bf A0, A1, A2},
and {\bf A3W} on cost functions introduced by Ma, Trudinger, and Wang.
The connectivity of contact sets of optimal transport maps follows as a direct corollary.
Our approach is quite different from the previous result of Loeper which he obtained as the
first step toward his H\"older regularity theory of potential functions, and which was based
upon approximation using the regularity theory of Ma, Trudinger, and Wang.
The result in this paper improves his result,
by relaxing certain geometrical assumptions on the domains and targets;
it also completes his H\"older regularity theory of potential functions on the round sphere,
by making it self-contained.
\end{abstract}

\section{Introduction}\label{S:intro}
For a real valued function $u$ on $\R^n$, its subdifferential (or subgradient) $\partial u(x)$ at $x$ is defined locally as
\begin{align*}
\partial u (x) = \{ p \in T_x\mathbf{R}^n \ | \ u(z) \ge
u(x) +  p \cdot z + o(|x-z|)  \ \hbox{as $ z \to x$} \}.
\end{align*}
It is clear that each $\partial u (x_0)$ is a convex set and if $u$ is differentiable at $x_0$,
$\partial u(x_0) = \{\nabla u (x_0)\}$. If $u$ is a convex function then for $p \in \partial u(x_0)$ the hyperplane
$H_p (x) = x \cdot p + u(x_0)$ supports $u$ globally, i.e. $u \ge H_p$ on $\mathbf{R}^n$.
In this paper, motivated by optimal mass transport theory, we are interested in the generalization of this fact for generalized convex functions, which are called cost-convex (or simply $c$-convex) functions, defined as
the generalized Legendre transformations
with respect to the cost function $c=c(x,y)$ on some product space $\Omega \times \Lambda$.
These are called $c$-Legendre transformations in Section~\ref{S:def}.

Ordinary convexity is the special case $c(x,y)=- x \cdot y$, and Brenier \cite{B} and others \cite{RR}\cite{KS}\cite{CP}\cite{Ra}\cite{M-duke}
showed that convex functions give potential functions for the optimal mass transport problem.
Namely, given two mass distributions $\mu_+$, $\mu_-$,
i.e. two Borel probability measures on $\mathbf{R}^n$,
assuming the measure $\mu_+$ does not charge sets of lower dimension,
the unique mapping $T_0$ which minimizes the average cost functional
\begin{align*}
\mathcal{C}(T) = \int_{\mathbf{R}^n} |x-T(x)|^2 d\mu_+
\end{align*}
among all Borel measurable maps which push-forward $\mu_+$ onto $\mu_-$,
in the sense that $\mu_- (E) = \mu_+ (T^{-1} (E))$ for all Borel set $E \subset \mathbf{R}^n$,
is actually given as the gradient mapping of some convex function $\phi$, i.e.
$T_0 = \nabla \phi$.  Here $T_0$ and its uniqueness are defined in the $\mu_+$-almost everywhere sense.
Brenier's result was generalized to other cost functions on Euclidean space
\cite{C2} \cite{GM}, to the Riemannian distance squared cost on
manifolds \cite{M} and to the Heisenberg group \cite{AR},
by Caffarelli, Gangbo and McCann, McCann, and Ambrosio and Rigot respectively.
Sufficient conditions for the map to be smooth were investigated by Delano\"e (for $n=2$) \cite{D},
Caffarelli (for $n\ge 2$) \cite{Ca2}\cite{Ca4}, and by Urbas \cite{U}, for the Euclidean
distance squared, and then
by Cordero-Erausquin and Delano\"e for the Riemannian distance squared on the torus 
\cite{Cr} and nearly flat manifolds \cite{D2}, and finally 
for other cost functions in a series of papers by
Ma, Trudinger, Wang \cite{MTW}\cite{TW}, and Loeper \cite{L}.
We refer the reader to the book \cite{V} by Villani for an introduction to optimal mass
transportation theory and its applications in partial differential equations and geometry.

For a $c$-convex function $\phi$ on the domain $\Omega$ one may define the $c$-subdifferential
$\partial^c \phi(x_0)$ at $x_0$ to consist of the gradient vectors $\nabla M_y (x_0)$ of the
functions $M_y (x) = -c(x, y) + \hbox{const.}$ which globally support $u$ at $x_0$,
meaning $u(z) \ge M_y (z)$ for all $z \in \Omega$ plus the equality $u(x_0) = M_y(x_0)$.
The mapping from $x_0$ to the contact set, namely the set of such points $y \in \Lambda$,
is called the $c$-contact mapping of $\phi$, and generalizes Brenier's subgradient mapping
for ordinary convex functions. Unlike $\partial \phi$, $\partial^c \phi$ is a global notion
depending on the whole domain of $c$, $\Omega \times \Lambda$, but for $\phi$ convex and $c(x,y) = -x\cdot y$, $\partial^c \phi =\partial \phi$ and one may naturally ask whether this indeed holds for other more general cost functions.
\begin{question}\label{Q:main}
For which cost functions $c$ on  domain $\Omega \times \Lambda$, is there the equivalence between the $c$-subdifferential $\partial^c \phi$ and the ordinary differential $\partial \phi$ for $c$-convex functions $\phi$?
\end{question}

An answer was given by Loeper \cite{L}: there is such an equivalence
$\partial^c \phi = \partial \phi$ under certain conditions on the function $c$,
introduced by Ma, Trudinger, and Wang \cite{MTW}\cite{TW}, who called them
{\bf A0, A1, A2}, and {\bf A3W} , and for bounded domains
$\Omega, \Lambda \subset \mathbf{R}^n$ satisfying certain geometric restrictions with respect
to $c$, $c^*$, where $c^*(x,y) =c(y,x)$,  namely that $\Omega$ and $\Lambda$ are mutually
uniformly strictly $c$, $c^*$-convex and $c$, $c^*$-bounded (see Section~\ref{S:Previous}).
Loeper's result is based upon an approximation argument which uses Trudinger and Wang's
$C^3$-regularity theory of potential functions of optimal transport for general cost
functions \cite{TW}, and those conditions on $c$ and on $\Omega$, $\Lambda$ which Loeper
needed were actually the ones originally imposed by Trudinger and Wang when they applied a continuity
method to the Monge-Amp\`ere type equations satisfied by the potential functions.
Loeper's result on $c$-convex functions was actually the first step toward his H\"older
regularity theory of potential functions for optimal transport with rough mass distributions,
which need not be absolutely continuous with respect to the Lebesgue measure,  under a strict
version {\bf A3S} of the condition {\bf A3W} (see \cite{L}). He also extended his H\"older
regularity theory to
the standard $n$-dimensional sphere with
cost given by its Riemannian distance squared, and to the reflector antenna problem, which
helped to motivate Wang's interest in the theory of optimal mass transportation \cite{W1}\cite{W2}
and has also been studied by many
others, e.g.\ \cite{CGH}\cite{KO}.

The sufficient condition {\bf A3W} of Ma, Trudinger, and Wang for regularity
was shown by Loeper \cite{L} to be
necessary for the continuity of optimal maps with smooth integrable mass
distributions whose Radon-Nikodym derivatives are bounded away from $0$ and $+ \infty$.
It is also necessary for the equivalence $\partial^c \phi = \partial \phi$ for $c$-convex
functions.

The goal of this paper is to provide a direct and geometric proof of the equivalence
$\partial^c \phi = \partial \phi$ between the cost-subdifferential and the ordinary
subdifferential for cost-convex functions $\phi$, under Ma, Trudinger, and Wang's conditions
{\bf A0, A1, A2}, and {\bf A3W}. Our result only assumes the mutual (not necessarily strict)
$c$ and $c^*$-convexity of $\Omega$ and $\Lambda$. For example we do not assume them
to be bounded,  and this allows us to deal with such domains as $\mathbf{R}^n$.
In fact we even assume slightly less (see Remark~\ref{R:weaker convex}), and our result
can easily be applied to domains in the round sphere $S^n$, including $S^n$ itself
(see Example~\ref{E:sphere} and Example~\ref{E:antenna}), thus it completes Loeper's
H\"older regularity theory of potential functions for optimal mass transport maps on
domains in $\mathbf{R}^n$ or in $S^n$ without relying on Trudinger and Wang's regularity result
or any other method.   As a consequence, it strengthens the results of
Cafferelli, Gutierrez, and Huang \cite{CGH}, by weakening their assumptions on the incoming
intensity in the reflector antenna problem.

In this research announcement,  we adopt the point of view that our domain $\Omega$ and
target $\Lambda$ are either subdomains of Euclidean space,  or else manifolds equipped
with Riemannian metrics which allow us to identify their tangent and cotangent spaces.
Unless otherwise stated, computations of higher derivatives are always carried out with
respect to Riemannian exponential coordinates centered at the point where the derivatives
are evaluated.  In a subsequent paper,  we abandon this
artifice and give a fully covariant development of the results for general
cost functions on smooth manifolds \cite{KM}.  Aside from providing conceptual clarity and
geometric insight, this point of view offers the possibility of extending the regularity
theory of Loeper (and presumably that of Ma, Trudinger and Wang) to other geometries.

The present manuscript is organized as the following:
In Section~\ref{S:def} we list the basic definitions and assumptions we need in this paper;
in Section~\ref{S:main} we give the precise statements of main results;
in Section~\ref{S:Previous} we summarize some of the relevant prior results by
Ma, Trudinger, Wang, and Loeper; throughout Section~\ref{S:Proof T:c-subdifferential} and
Section~\ref{S:Proof T:monotonicity} we prove our main results; in Section~\ref{S:examples}
we discuss some examples of interest.

\subsection*{Acknowledgment}
We thank Gregoire Loeper, Neil Trudinger, and Xu-Jia Wang for their helpful discussions and
comments, and their recent preprints.  In particular, while writing this article,
we learned from Trudinger and Wang that they had independently relaxed some of the restrictions on
their earlier theory using a different method \cite{TW2}.
There they indicated the desirability of a direct approach to such questions,
which we believe the present manuscript provides.  We are grateful to
Adrian Nachman and the 2006-07 participants of
Fields Analysis Working Group, for the stimulating environment which they helped
to create.
This research was supported in part by Natural Sciences and Engineering
Research Council of Canada Grant 217006-03 and United States National Science Foundation
Grant DMS-0354729.


\section{Definitions }\label{S:def}
In this section we list some (modified) definitions introduced by  Ma, Trudinger, Wang, and
Loeper \cite{MTW}\cite{TW}\cite{L}. In many cases, we follow the notation of \cite{L},
despite the fact that it is at odds with the notation of \cite{GM}.
Let $\Omega$, $\Lambda$ be two (not necessarly bounded) domains in Riemannian manifolds
$M$ and $N$ respectively; for example $M=N=\mathbf{R}^n$ with the Euclidean metric.
\begin{definition}
A cost function $c : \Omega \times \Lambda \to \mathbf{R}$ is said to satisfy\\
{\bf(A0)}  if $c \in C_{\text{loc}}^4 (\overline{\Omega} \times \overline{\Lambda})$,\\
{\bf (A1)} if for all $x \in \Omega$  the mapping
$ y \to -\nabla_x c(x,y)$ is injective on $\overline{\Lambda}$ and  for all $y \in \Lambda$ the mapping $x \to -\nabla_y c(x,y)$ is injective on $\overline{\Omega}$,\\
{\bf (A2)} if $\det D^2_{xy} c(x,y) \ne 0$, for all $x \in \overline{\Omega}$ and $y  \in \overline{\Lambda}$ .
\end{definition}

\begin{definition}\label{c-exp}{\bf($c$-Exponential Maps)} Under the assumption {\bf A1}, we define the $c$ and $c^*$-exponential maps and their domains by
\begin{align*}
&c\text{-Exp}_x \ p = y \hbox{\ if $p = -\nabla_x c(x, y)$,}\\
&\text{dom}(c\text{-Exp}_x ) =\{ p \in T_x\Omega \  | \  p = -\nabla_x c(x, y) \hbox{ \ \ for some $y \in \Lambda$}  \}.
\end{align*}
and
\begin{align*}
&c^*\text{-Exp}_y \ q = x \hbox{\ if $q = -\nabla_y c(x, y)$,}\\
&\text{dom}(c^*\text{-Exp}_x ) =\{ q \in T_y \Lambda \ | \ q= -\nabla_y c(x, y) \hbox{ \ \ for some $x \in \Omega$} \}.
\end{align*}
The assumption {\bf A2} implies that $c$ and $c^*$-exponential maps are local diffeomorphisms.
\end{definition}

Throughout this paper we will assume {\bf A0, A1, A2} for the cost function
$c : \Omega \times \Lambda \to \mathbf{R}$ and we define $c^*$ by  $c^*(y,x) =c(x,y)$.
The following simple observation will be useful.

\begin{lemma}\label{L:symmetry}{\bf (Symmetry)}
For $\eta \in T_x \Omega$, $\xi \in T_y \Lambda$,
\begin{align*}
(D^2 _{x, y} c (x, y))^{-1} \eta \cdot \xi =
\eta \cdot (D^2 _{y, x} c (x,y) )^{-1}\xi \end{align*}
i.e.
\begin{align*}
 [d \ c\hbox{-Exp}_x ]_{(c\hbox{-Exp}_x)^{-1}(y) }(\eta) \cdot \xi = \eta \cdot
 [d \ c^*\hbox{-Exp}_y]_{(c^*\hbox{-Exp}_y)^{-1} (x)} (\xi)
\end{align*}
where `$\cdot$' means the Riemannian inner product.
\end{lemma}
\begin{proof}
The proof is straightforward.
\end{proof}

\begin{definition}\label{D:cost sec A3}{\bf ($c$-sectional curvature and A3W, A3S)}
Define the \emph{$c$-sectional curvature} for $(x,y) \in \Omega \times \Lambda$, $\eta, \xi \in T_x \Omega$,
by
\begin{align*}
 \mathfrak{S}_c (x, y) (\eta, \xi) :
 = \frac{d^2}{dt^2}\Big{|}_{t=0} [-D^2_{xx} c] (x, c\hbox{-Exp}_x (p + t \xi)) \ \eta \  \eta
\end{align*}
for $p=-\nabla_x c(x,y) = [c\hbox{-Exp}_x ]^{-1}(y)$.
We say that $c$ satisfies {\bf A3W} (resp. {\bf A3S}) if  there exists a constant $C_0 \ge 0$ (resp. $>0$) such that for all  $(x,y) \in \Omega \times \Lambda$,
\begin{equation}\label{c-sectional curvature}
   \mathfrak{S}_c (x,y) (\eta, \xi) \ge C_0 |\eta|^2 |\xi|^2   \ \  \hbox{for all $\eta \perp \xi$.}
\end{equation}

\end{definition}

\begin{definition}\label{D:c-convexity}
{\bf ($c$-convexity)}
We say $\Lambda$ is \emph{$c$-convex with respect to $x \in \Omega$}
if $ (c\text{-Exp}_x )^{-1}\Lambda$ is convex as a set in the Euclidean space $T_x \Omega$.
 We say $\Lambda$ is \emph{$c$-convex with respect to $\Omega$}
 if it is $c$-convex with respect to all $x \in \Omega$.
We define similarly the $c^*$-convexity of $\Omega$.
\end{definition}

\begin{definition}\label{D:c-convex function}
{\bf ($c$-convex function indexed by a set)}
We define the \emph{c-Legendre transformation} $L^c v$ of $v:\Lambda \to \mathbf{R}\cup \{+\infty \}$ to
be the function on $\Omega$ given by
\begin{align*}
L^c (v) (x)= \sup_{y \in \Lambda} -c(x,y) - v(y).
\end{align*}
We define a function $\phi$ on $\Omega$ to be a \emph{$c$-convex function indexed by  $\Lambda$}, if
\begin{align} &\hbox{$\phi= L^c(v)$ for some function $v$ on $\Lambda$}\\\nonumber
& \hbox{and}\\\label{additional condition}
&\hbox{for each $x \in \Omega$ there is $y \in \overline{\Lambda}$ such that $\phi(x) + v(y) = -c(x,y)$.}
\end{align}

Setting $c^* (x,y) = c(y,x)$ we may define
the $c^*$-Legendre transformation $L^{c^*}$ from functions on $\Omega$ to the functions on $\Lambda$, and  $c^*$-convex functions indexed by $\Omega$, similarly.
\end{definition}

\begin{remark}
\eqref{additional condition} is not explicitly assumed for $c$-convexity in \cite{MTW}\cite{TW}\cite{L} since the domains there $\Omega$ and $\Lambda$ are assumed to be bounded, and hence \eqref{additional condition} is automatically satisfied due to compactness of $\overline{\Omega} \times \overline{\Lambda}$ and continuity of $c$.
\end{remark}


 \begin{definition}{\bf ($c$-contact mapping)}
Let $u$ be  a $c$-convex function on $\Omega$ indexed by $\Lambda$. We may continuously extend it to the closure $\overline{\Omega}$ of $\Omega$.  We define the set-valued mapping
$G_{u}$ from $\overline{\Omega}$ to $\overline{\Lambda}$ by
\begin{align*}
G_{u}(x):= \{ y \in \overline{\Lambda} \ | \ \  u (x) + L^{c^*}(u) (y) = -c(x, y) \}
\hbox{  for each $x \in \overline{\Omega}$}.
\end{align*}
$G_u$ is called \emph{the $c$-contact mapping} (or simply \emph{the contact mapping}) of $u$, and the set $G_u (x)$ is called \emph{the contact set} of $u$ at $x$.
We can also define in a similar way the contact map $G_v$ from $\overline{\Lambda}$ to $\overline{\Omega}$, for a $c^*$-convex function $v$ on $\Lambda$ indexed by $\Omega$. Then, it is clear that for any  $x \in \overline{\Omega}, \ y \in \overline{\Lambda}$ if $v=L^{c^*}(u)$,
\begin{align*}
G_{u}(x) = (G_{v})^{-1} (x), \ \ \
G_{v}(y) = (G_{u})^{-1} (y) .
\end{align*}
\end{definition}


\begin{definition}{\bf (subdifferential and $c$-subdifferential)}
Assume that  $\Omega$ is given a Riemannian inner product $\cdot$
and the distance function $d$ induced by this metric.
Let $u$ be a $\mathbf{R}\cup\{+\infty\}$-valued function on $\Omega$. The \emph{subdifferential} $\partial u$ is defined locally for each $x \in \Omega$  by
\begin{align*}
\partial u (x) = \{ p \in T_x \Omega \ | \ u(z) \ge
u(x) +  p \cdot \exp_x^{-1} z + o(\text{d}(x,z)) \hbox{  as $z \to x$} \}.
\end{align*}
Note that the subdifferential is always a convex set though it may possibly be empty. However,
 for a $c$-convex function $u$ indexed by $\Lambda$, this set is nonempty, and we may further define the \emph{$c$-subdifferential} $\partial^{c} u $  for each $x_0 \in \Omega$ by
\begin{align*}
\partial^{c} u (x_0) :=
\big{\{} p \in T_{x_0}\Omega \ | \ p =- \nabla_{x=x_0} c (x, y), \ y \in G_{u} (x_0) \big{\}}.
 \end{align*}
 One should note that $\partial^c u (x) \subset \partial u (x)$ for any $x \in \Omega$ and that $\partial^c u$ depends globally on the sets $\Omega$, $\Lambda$.
 \end{definition}

\section{Main Results}\label{S:main}
Using the definitions given in the previous section, we can now precisely state our main result.
It was formulated by Loeper \cite{L}, who proved it under various technical hypotheses
described in the section to follow using a result of Trudinger and Wang \cite{TW}.

\begin{theorem}\label{T:c-subdifferential} {\bf (A3W $\Leftrightarrow$ CSIS)} Let $\Omega$, $\Lambda$ be two (not necessarly bounded) domains in Riemannian manifolds $M$ and $N$, resp. (e.g. $M=N=\mathbf{R}^n$, with Euclidean metric).
Assume that the cost $c: \Omega \times \Lambda \to {\bf R}$ satisfies {\bf A0, A1, A2} and that
 $\Lambda$ is $c$-convex with respect to $\Omega$, and $\Omega$ is $c^*$-convex with respect to $\Lambda$. Then the following are equivalent.

1. {\bf (A3W)} $c$ satisfies {\bf A3W}.

2. {\bf (CSIS: c-subdifferential is subdifferential)}
For any  $c$-convex function $\phi$ on $\Omega$ indexed by $\Lambda$,
\begin{align*}
\partial^c \phi = \partial \phi .
\end{align*}

\end{theorem}

From this theorem, we can deduce the connectivity of contact sets of $c$-convex functions when we assume {\bf A3W}.
\begin{corollary}\label{C:connectivity} {\bf (connectivity of contacts sets)}
Under the assumption of Theorem~\ref{T:c-subdifferential}, if further {\bf A3W} holds, then
for any $c$-convex functions $\phi$ on $\Omega$ indexed by $\Lambda$, the contact set
$G_\phi (x) \subset \overline{\Lambda}$ of $\phi$ is $c$-convex (and so connected) for all $x \in \Omega$.
\end{corollary}

The following corollary is useful, especially when $\Omega =\mathbf{R}^n$ (see Example~\ref{E:convex f g} and Example~\ref{E:log distance}), and it is parallel with Caffarelli's result \cite{Ca2} for Euclidean distance squared cost, which required the convexity only of the target domain $\Lambda$.

\begin{corollary}\label{C:subdomains}{\bf (arbitrary subdomains)}
Put the same assumption on $c$, $\Omega$, $\Lambda$ as in Theorem~\ref{T:c-subdifferential}. Then for any arbitrary subdomain $\Omega_1 \subset \Omega$ and for a subset $\Lambda_1 \subset \Lambda$ $c$-convex with respect to $\Omega_1$, and for any $c$-convex function $\phi$ on $\Omega_1$ indexed by $\Lambda_1$,  we have $\partial^c \phi = \partial \phi$.
\end{corollary}
\begin{proof}
The proof is given at the end of Section~\ref{S:Proof T:c-subdifferential}.
\end{proof}

\section{Previous results: Ma, Trudinger, Wang, and Loeper}\label{S:Previous}
We now discuss the previous results of Ma, Trudinger, Wang, and Loeper.
The following conditions were critical in their initial studies \cite{MTW}\cite{TW}\cite{L}.
Our approach avoids many of these restrictions,
so that our arguments in Section~\ref{S:Proof T:c-subdifferential} and Section~\ref{S:Proof T:monotonicity}
only assume (not necessarily uniformly strict) $c$, $c^*$-convexity --- an improvement
which, like Corollary \ref{C:subdomains}, was
attained independently in the parallel work of Trudinger and Wang \cite{TW2},

We denote
\begin{align}\label{tensor A}
\mathcal{A}(x,p) = -D^2_{xx} c (x,y) \ \hbox{ for $y$ such that $-\nabla_x c(x,y)= p$}.
\end{align}
\begin{definition}\label{uniformly c-convex} {\bf (Uniform strict $c$, $c^*$-convexity) }We say that  $\Omega$ is \emph{uniformly strictly $c^*$-convex} with respect to $\Lambda$ if its inverse $c^*$-exponential image
$[c^*\text{-Exp}_y ]^{-1} (\Omega)$ is uniformly strictly convex with respect to $y \in \Lambda$, i.e., if $\Omega$ is $c^*$-convex with respect to $\Lambda$, $\partial \Omega \in C^2$ and there exists a constant $\delta_0 >0$  such that
\begin{align*}
 [D_{i} \gamma_j (x) - D_{p_k}
 \mathcal{A}_{ij} (x, p) \gamma_k ] \tau_i \tau_j \ge \delta_0
\end{align*}
for all $x\in \partial\Omega$, $p\in -\nabla_x c(x, \Lambda) \subset T_x \Omega$, unit tangent vector $\tau$ and outer unit normal vector $\gamma$. Here $\mathcal{A}$ is given in \eqref{tensor A}.
We similarly define the uniform strict $c$-convexity of $\Lambda$ with respect to $\Omega$.
\end{definition}

\begin{definition}\label{c-boundedness}
{\bf ($c$, $c^*$-boundedness)}
We say that $\Omega$ is \emph{$c^*$-bounded}, with respect to $\Lambda$ if there exists some function $\varphi \in C^2 (\Omega)$, satisfying for some constant $\delta_1 >0$,
\begin{align*}
[D^2_{ij} \varphi (x) - D_{p_k} \mathcal{A}_{ij} (x,p) D_k \varphi ] \xi_i \xi_j \ge \delta_1 |\xi|^2
\end{align*}
for all $x \in \Omega$ and for all  $p\in -\nabla_x c(x, \Lambda) \subset T_x \Omega$. Here $\mathcal{A}$ is given in \eqref{tensor A}.
We similarly define the $c$-boundedness of $\Lambda$.
\end{definition}

\begin{remark*}
In \cite{TW}, Trudinger and Wang mean by $c$-convexity and $c$-boundedness of $\Omega$, the $c^*$-convexity and $c^*$-boundedness in this paper. We follow the convention given in \cite{L}.
\end{remark*}
\begin{remark*} (\cite{TW}, p19)
Any sufficiently small ball $\Lambda$ is $c$-bounded with respect to a bounded set $\Omega$. \end{remark*}

In \cite{MTW} \cite{TW}, Ma, Trudinger and Wang obtained the following regularity result for
$c$-convex functions as solutions to optimal transport problems.
\begin{theorem}\label{T:Trudinger-Wang} ({\bf Ma, Trudinger, and Wang.} \cite{MTW} \cite{TW} see also \cite{L} ,p13)
Let $\Omega$, $\Lambda$ be two bounded domains in $\mathbf{R}^n$. Let $c$ satisfy {\bf A0,A1,A2}, and {\bf A3W} on $\Omega \times \Lambda$. Assume that $\Omega$, $\Lambda$ are strictly uniformly $c, c^*$-convex with respect to each other, either $\Lambda$ is $c$-bounded with respect to $\Omega$ or $\Omega$ is $c^*$-bounded with respect to $\Lambda$.  Then for two probability measures
$\mu_0$, $\mu_1$  on $\Omega$, $\Lambda$ having densities $\rho_1 \in C^2 (\overline{\Omega})$, $\rho_2 \in C^2(\overline{\Lambda})$ which are bounded away from $0$, and for any  $c$-convex function $\phi$ on $\Omega$, with $G_\phi (\Omega) \subset \overline{\Lambda}$, $(G_\phi)_\# \mu_0 = \mu_1$, we have $\phi \in C^3 (\overline{\Omega})$.
\end{theorem}

Using this regularity result, Loeper showed  the equivalence
between $c$-subdifferentials and ordinary subdifferentials
(under the same assumptions of Trudinger and Wang), by approximating $c$-convex functions by smooth ones.
\begin{theorem}\label{T:Loeper c-subdifferential}
({\bf Loeper} \cite{L},p13, Theorem 3.1.)
Under the same assumption on the cost $c$, domains $\Omega$, $\Lambda$ as in Theorem~\ref{T:Trudinger-Wang} (especially with $c$ or $c^*$-boundedness),
we have for any $c$-convex function $\phi$ on $\Omega$, with $G_\phi (\Omega) \subset \overline{\Lambda}$,
that
\begin{align*}
\partial^c \phi (x) = \partial \phi  (x) \ \ \hbox{\ for all $x \in \Omega$}, \end{align*}
and so the set $\partial^c \phi(x) \in T_x\Omega$ is convex for all $x \in \Omega$.
\end{theorem}

\section{Proof of Theorem~\ref{T:c-subdifferential}}\label{S:Proof T:c-subdifferential}
In this section, we illustrate a blueprint of the proof of Theorem~\ref{T:c-subdifferential}.

We first start with the geometric setting of the proof. Pick $x_m \in \Omega$ and two different points $y_0, y_1 \in \Lambda$ and their inverse $c$-exponential images $p_0, p_1 \in T_{x_m} \Omega$ , $p_0 \ne p_1$.
Let $p_\theta = (1-\theta) p_0 + \theta p_1$, $\theta \in \mathbf{R}$, be the line in $T_{x_m} \Omega$ passing through $p_0$ with direction $p_1 -p_0$.
For such $p_\theta \in \text{dom}( c\text{-Exp}_{x_m})$, consider its $c$-exponential image $y_\theta$, i.e.
\begin{align}\label{p_theta}
 p_\theta = -\nabla_x c(x_m, y_\theta)  \in T_{x_m} \Omega, \ \ \   y_\theta = c\text{-Exp}_{x_m}(p_\theta).
\end{align}
Then, by the condition {\bf A2}, we see that $\frac{d}{d\theta}y_\theta \ne 0$.
If such $y_\theta$ is defined at least for $\theta \in [0,1]$, then we denote by $[y_0, y_1]_{x_m}$ the curve $\{ y_\theta \}_{0 \le \theta \le 1}$ and call it  the \emph{$c$-segment} with respect to $x_m$ joining $y_0$ to $y_1$ \cite{L}.
Such $[y_0, y_1]_{x_m}$ exists, for example if $\Lambda$ is $c$-convex with respect to $x_m$.

\begin{definition}\label{D:sliding mountain}{\bf (Sliding and double mountains)}
Given $x_m \in \Omega$ and $y_0,y_1 \in \Lambda$, define the function $f_\theta$ on $\Omega$,
which we call the \emph{sliding mountain} between $y_0$,$y_1$ centered at $y_\theta$ and normalized at $x_m$,
by
\begin{align}\label{f_theta}
f_\theta(x) =-c(x, y_\theta) + c(x_m , y_\theta)
\end{align}
where $y_\theta$ is the point corresponding to $\theta \in [0, 1]$ in the $c$-segment $[y_0, y_1 ]_{x_m}$.
We call $\max[f_0 , f_1]$ the \emph{double mountain}  of $y_0$, $y_1$, normalized at $x_m$.
\end{definition}

\begin{definition}\label{D:S theta}{\bf (The set $S_\theta$ and $S^+ _\theta$)} We define  for the sliding mountain $f_\theta$,
\begin{align*}
\hbox{the level set $S_\theta: =\{x \in \Omega \ | \ \frac{d}{d\theta} f_\theta (x) = 0 \}$  and }\\
\hbox{the super-level set
$S^+_\theta:= \{ x \in \Omega \ | \ \frac{d}{d\theta} f_\theta (x) \ge 0 \}$} .
\end{align*}
\end{definition}

In this notation we have the following key observation whose proof is straightforward.
\begin{lemma}\label{L:DMASM=MONO}{\bf (DASM $\Leftrightarrow$ Monotonicity)}
For $x_m \in \Omega$, suppose that
$y_0, y_1 \in \Lambda$ is joined by the $c$-segment $[y_0 , y_1]_{x_m}$ with respect to $x_m$.
Then (1) is equivalent to (2),  and (1')\ is equivalent to (2'):

(1) {\bf (DASM: Double Mountain Above Sliding Mountains)} \begin{align}\label{double-sliding}
f_\theta (x) \le \max[f_0 , f_1](x) ,  \ \hbox{  for all $(\theta, x) \in ]0,1[ \times \Omega$} .
\end{align}

(2) {\bf (Monotonicity)} The
super-level set
$S^+_\theta$ monotonically increases in $\theta$, i.e.
\begin{align}\label{monotonicity}
 S^+_{\theta_1} \subset S^+_{\theta_2}, \ \ \hbox{\  for any  $0 \le \theta_1 \le \theta_2 \le 1$}.
\end{align}

(1') {\bf (Strict DASM) } (\ref{double-sliding}) holds, and the inequality is strict unless $x=x_m$.

(2') {\bf (Strict Monotonicity)} (\ref{monotonicity}) holds and
$S_{\theta_1} \cap S_{\theta_2} = \{ x_m \}$ if $\theta_1 < \theta_2$.

\end{lemma}
\begin{remark}
 In this lemma, we do not assume $c$, $c^*$-convexity on $\Omega$, $\Lambda$.
\end{remark}

  \begin{remark}\label{R:double-sliding}{\bf (DASM $\Leftrightarrow$ CSIS)}
 Note that \eqref{double-sliding} is equivalent to the convexity of the $c$-subdifferential $\partial^c \max[f_0, f_1] (x_m)$ at $x_m$, and so to $\partial^c \max[f_0, f_1] (x_m)=\partial \max[f_0, f_1] (x_m)$.  In fact, it is easy to see as observed by Loeper (see Proposition 2.12 in \cite{L}) that  for $\Lambda$ being $c$-convex with respect to $\Omega$, the property
\eqref{double-sliding} holds for every $x_m, x \in \Omega$, and for all $c$-segments $[y_0, y_1]_{x_m} \subset \Lambda$ if and only if
 $\partial^c \phi = \partial \phi$ for all $c$-convex functions $\phi$ on $\Omega$ indexed by $\Lambda$.
\end{remark}
\begin{remark}\label{R:Loeper}{\bf (DASM $\Rightarrow$ A3W)}
 Loeper used Taylor expansion to prove the important result that
 \begin{align*}
 \hbox{ \eqref{double-sliding} holds locally near $x_m$} \Rightarrow {\bf A3W}  .
\end{align*}
(See \cite{L} pp20--21.)
\end{remark}

From Lemma~\ref{L:DMASM=MONO}, Remark~\ref{R:double-sliding}, and Remark~\ref{R:Loeper}, it is clear that to prove  Theorem~\ref{T:c-subdifferential},
it suffices to prove either {\bf (DASM)} or {\bf (Monotonicity)} for all $x_m$, $x$, $[y_0 , y_1]_{x_m}$ while assuming {\bf A3W}.
We will do this by proving the following theorem in Section~\ref{S:Proof T:monotonicity}.

\begin{theorem}\label{T:monotonicity}{\bf (A3 $ \Rightarrow$ Monotonicity)}
Under the same assumption on $c$, $\Omega$, $\Lambda$ as in Theorem~\ref{T:c-subdifferential},
if $c$ satisfies {\bf A3W} then {\bf(Monotonicity)} (\ref{monotonicity})
holds for all $x_m \in \Omega$ and $y_0 , y_1 \in \Lambda$.  Similarly {\bf A3S} implies that
{\bf (Strict Monotonicity)} holds for all $x_m \in \Omega$ and $y_0 , y_1 \in \Lambda$.
\end{theorem}

We emphasize that our proof of Theorem~\ref{T:monotonicity}, thus the direction
$(1) \Rightarrow (2)$ in Theorem~\ref{T:c-subdifferential}, is elementary and geometric.
It is not based on the previous results of Ma, Trudinger, Wang, and Loeper.

We finish this section by giving the proof of Corollary~\ref{C:subdomains}.

\subsubsection*{Proof of Corollary~\ref{C:subdomains}}
It is clear that the property {\bf (DASM)} \eqref{double-sliding} has monotonicity with respect to the domain: For $\Omega_1 \subset \Omega$, $\Lambda_1 \subset \Lambda$ with $\Lambda$, $\Lambda_1$ being $c$-convex with respect to $\Omega$, $\Omega_1$, respectively,
if \eqref{double-sliding} holds for every $x_m, x \in \Omega$, and for all $c$-segments $[y_0 , y_1 ]_{x_m} \subset \Lambda$, then the same holds for the pair $\Omega_1, \Lambda_1$.
Therefore we can conclude using Remark~\ref{R:double-sliding} and Theorem~\ref{T:c-subdifferential} that
$\partial^c \phi = \partial \phi$ for all $c$-convex functions $\phi$ on $\Omega_1$ indexed by $\Lambda_1$.

\section{Proof of Theorem~\ref{T:monotonicity} {\bf(A3 $\Rightarrow$ Monotonicity)}}\label{S:Proof T:monotonicity}
We divide the proof into three steps.
We will first set up the necessary notation and show some preliminary geometric facts in the subsection~\ref{SS:geometry}.
Then, we will show in the subsection~\ref{SS:rate of expansion} that the rate of expansion of $S^+_\theta$
 along $S_\theta$ (see Definition~\ref{D:sliding mountain} and Definition~\ref{D:S theta}) has the same sign as the quantity
$\frac{d^2}{d\theta^2}f_\theta$ on $S_\theta$. In the subsection~\ref{SS:positivity} we will show the
non-negativity of this expansion rate to conclude the desired monotone dependence of the rising region
$S^+_\theta$ on $\theta \in [0,1]$ (see Proposition~\ref{P:double f positive}),
and it will finish our proof of Theorem~\ref{T:monotonicity}.

\subsection{Geometric preliminaries}\label{SS:geometry}
We first recall the definitions of $p_\theta, y_\theta$ from Section~\ref{S:Proof T:c-subdifferential}.
For each $y_\theta$, we define
\begin{align}\label{q_theta}
 q_\theta = -\nabla_y c(x_m , y_\theta)  \in T_{y_\theta} \Lambda , \ \ \ \ x_m = c^*\text{-Exp}_{y_\theta} (q_\theta)
\end{align}
We denote $\dot{y}_\theta = \frac{d}{d\theta} y_\theta$,
$\ddot{y}_{\theta} = \frac{d^2}{d\theta^2} y_\theta$, and we also use $\cdot$ to denote
the Riemannian inner product. Our calculations will be performed in geodesic normal coordinates of $x \in \Omega$ and $y_\theta \in \Lambda$.

Now, let's discuss the geometric properties of $S_\theta$.
We first note that
\begin{align*}
D_x \frac{d}{d\theta} f_\theta \cdot \xi
= \big{(}-D^2_{y,x} c(x, y_\theta) \ \xi \big{)} \cdot \dot{y}_\theta,
\ \ \hbox{for any tangent vector $\xi$ at $x$},
\end{align*}
and by {\bf A2} that  $D_x \frac{d}{d\theta} f_\theta$ is nowhere vanishing on $\Omega$. Thus we see the zero level set $S_\theta$ of $\frac{d}{d\theta} f$  is a smooth $(n-1)$-dimensional submanifold in the $n$-dimensional domain $\Omega$. It is clear that $x_m \in S_\theta$ and that
\begin{align}\label{S_theta}
S_\theta = c^*\text{-Exp}_{y_\theta} ((q_\theta + W_\theta)  \cap \text{dom}(c^*\text{-Exp}_{y_\theta})),
\end{align}
for the $(n-1)$-dimensional subspace
\begin{align}\label{W theta}
W_\theta = \{ w \in T_{y_\theta} \Lambda \ |
\  w \cdot  \dot{y}_\theta =0 \} \ \ \subset T_{y_\theta} \Lambda     .
\end{align}
We parametrize by $w$ the points in $W_\theta \subset T_{y_\theta} \Lambda$
and for each $w$ let $x_w $ be the point in $S_\theta  \subset \Omega$ given by
\begin{align}\label{parameter w}
w = -D_y c(x_w  , y_\theta) + q_\theta .
\end{align}
At each point $x_w \in S_\theta$, denote by $p_\theta (x_w)$ the inverse $c$-exponential image of $y_\theta$ given by
 \begin{align}\label{p theta}
p_\theta (x_w ) = -D_x c(x_w, y_\theta),
\end{align}
and denote
\begin{align}\label{dot p theta}
\dot{p}_\theta (x_w) = \frac{d}{d\theta} p_\theta (x_w).
\end{align}
Note that $w \cdot \dot{y}_\theta =0$, and by Lemma~\ref{L:symmetry}
$\dot{p}_\theta (x_w)$  is orthogonal to $S_\theta $ at each $x_w$, i.e.
 \begin{align}\label{orthogonal to S_theta}
 \dot{p}_\theta (x_w) \cdot D_w x_w =0 ,
\end{align}
when the derivative $D_w x_w$ is taken for any direction of $w$.
In particular, we see that $S_\theta$ is orthogonal to $p_1 - p_0$ at $x_m$ for all $\theta$.

As our key ingredients,
we consider for each $\theta$ the functions (on $\Omega$)
\begin{align*}
 \frac{d}{d\theta} f_\theta (x)
&=\big{(} -D_y c(x, y_\theta) + D_y c(x_m, y_\theta) \big{)} \ \dot{y}_\theta , \\
\frac{d^2}{d\theta^2} f_\theta (x) & =
\big{(} -D^2_{yy} c (x, y_\theta) + D^2_{yy} c(x_m , y_\theta) \big{)} \dot{y_\theta} \ \dot{y_\theta}   \\ & \ + \big{(} -D_y c(x,y_
\theta) + D_y c (x_m, y_\theta) \big{)} \ \ddot{y}_{\theta}.
\end{align*}
Note that at $x_m$,
$f_\theta (x_m), \frac{d}{d\theta}f_\theta (x_m), \frac{d^2}{d\theta^2} f_\theta (x_m ) = 0$,
and for all $\xi \in T_{x_m} \Omega$,
\begin{align*}
 D_x \frac{d^2}{d\theta^2} f_\theta (x_m) \xi  = &
 -D_x D^2_{yy} c (x_m, y_\theta) \  \xi \ \dot{y_\theta}  \dot{y_\theta}   \\ & \
 -D^2_{yx} c(x_m,y_
\theta)  \  \xi  \   \ddot{y}_{\theta}.
\end{align*}
And, by differentiating $(1-\theta)p_0 + \theta p_1  = -D_x c(x_m , y_\theta)$ twice with respect to $\theta$, we see
\begin{align}\label{zero y_theta}
0= -D^2_{yy} D_x c(x_m, y_\theta) \dot{y}_\theta \dot{y}_\theta
- D^2_{xy}  c(x_m, y_\theta)  \ddot{y}_\theta
\end{align}
and that
\begin{align}\label{gradient double f_theta at x_m}    D_x \frac{d^2}{d\theta^2} f_\theta (x_m) =0 .
\end{align}

\subsection{The rate of expansion of $S^+_\theta$ along $S_\theta$ has the same sign as the quantity
$\frac{d^2}{d\theta^2}f_\theta$ on $S_\theta$}\label{SS:rate of expansion}
This result comes form the following lemma, by letting $g_t = \frac{d}{d\theta}f_\theta \Big{|}_{\theta =t}$.
\begin{lemma}
Let $g : (-\epsilon, \epsilon) \times U \to \mathbf{R}$ be a $C^2$ function,
where $U$ is an open domain in $\mathbf{R}^n$.
Denote $g_t(x)= g(t, x)$ and
suppose the gradient $\nabla_x g_t$ is non-vanishing on $U$ for $t \in (-\epsilon, \epsilon)$.
Then the moving level set $L_t := \{ x \in U \ | \ g_t (x) = 0 \}$ is locally parametrized
for small $t$ by $X(t,x)$, $x \in L_0$, where
we define $X(t, x)$ by solving the following first order ODE :
\begin{align*}
\left\{%
 \begin{array}{ll}
\frac{d}{dt}X(t,x)= -\frac{\dot{g}_t(X(t,x))}{\|\nabla_x g_t (X(t,x)) \|^2} \nabla_x g_t (X(t,x)) & \ \ \\[1ex]
  X(0, x) = x  ,
      \end{array}%
  \right.
\end{align*}
where we denote $\dot{g}_t := \frac{\partial}{\partial t} g_t$.
Moreover, the expansion rate of the super level set $L^+ _t := \{ x \in U \ | \ g_t (x) \ge 0 \}$ along $L_t$ in this parametrization is given by
\begin{align}\label{expansion L+}
-\frac{\nabla_x g_t (x)}{\|\nabla_x g_t (x) \|} \cdot \frac{d}{dt}X(t,x) =  \frac{\dot{g}_t (x)}{ \|\nabla_x g_t (x) \|}
\end{align}

\end{lemma}
\begin{proof}
Since
\begin{align*}
 \frac{d}{dt}g_t (X(t,x)) =  \dot{g}_t (X(t,x)) +
 \nabla_x g_t (X(t, x)) \cdot \frac{d}{dt} X(t,x)
 =  0 ,
\end{align*}
it is clear that this $X(t, x)$ is indeed the desired parametrization of $L_t$.  \eqref{expansion L+} is clear from the construction of $X(t, x)$.
\end{proof}

\subsection{The positivity of $\frac{d^2}{d\theta^2}f_\theta \Big{|}_{S_\theta}$}\label{SS:positivity}
Now we will show that
$\frac{d^2}{d\theta^2}f_\theta \Big{|}_{S_\theta} \ge 0$ under {\bf A3W}. Under {\bf A3S} condition, the equality holds only at $x_m$, so it is the only common point of the expanding boundaries $S_\theta$. This result together with the result in the subsection~\ref{SS:rate of expansion} will finish the proof of Theorem~\ref{T:monotonicity}.

\begin{proposition}\label{P:double f positive}{\bf (Positivity of $\frac{d^2}{d\theta^2}f_\theta$ on $S_\theta$)}
Suppose that the cost function $c$ satisfies {\bf A0, A1}, \& {\bf A2}.
Assume further that $\Omega$ is $c^*$-convex with respect to $\Lambda$
and that the $c$-segment $\{y_\theta\}_{0 \le \theta \le 1}$ with respect to $x_m$ is in $\Lambda$.
Then we have:
\begin{align*} &(1). \hbox{\ \ \ under {\bf A3W}, $\frac{d^2}{d\theta^2} f_\theta \Big{|}_{S_\theta} \ge 0$ .}  \\
&(2). \hbox{\ \ \ under {\bf A3S}, $\frac{d^2}{d\theta^2} f_\theta \Big{|}_{S_\theta} (x) > 0$ for $x \ne x_m$.}
\end{align*}
\end{proposition}

To prove this proposition, let's consider the following key lemma whose proof is given at the end of this section.
\begin{lemma}\label{L:double double f_theta}
{\bf (Sliding mountain and $c$-curvature)}
We recall the $c$-sectional curvature $\mathfrak{S}_c$ in Definition~\ref{D:cost sec A3} and the parametrization $w$ of the set $W_\theta$ of tangent vectors in $T_{\dot{y}_\theta}\Lambda$ orthogonal to $\dot{y}_\theta$ as given in \eqref{W theta}. Let  $x_w$ be the points in $S_\theta \subset \Omega$ defined as in \eqref{parameter w}, and let $p_\theta$, $\dot{p}_\theta$ be the vectors defined in \eqref{p theta} and \eqref{dot p theta}. Fix a tangent vector $\eta \in T_w W_\theta$, and let
$D_w$ be the directional derivative in $w$ of the direction $\eta$ and let $D^2_{ww}$ be
the second iteration of $D_w$. Then
 \begin{align*}
 D^2_{ww} \frac{d^2}{d\theta^2} f_\theta (x_w)
= \mathfrak{S}_c (x_w , y_\theta) (D_w x_w , \dot{p}_\theta(x_w) ).
 \end{align*}
 \end{lemma}

\begin{proof}[Proof of Proposition~\ref{P:double f positive}]

Denote
\begin{align*}
V_\theta:=(q_\theta + W_\theta ) \cap \text{dom}(c^*\text{-Exp}_{y_\theta})
\end{align*}
and note that $S_\theta = c^*\text{-Exp}_{y_\theta} (V_\theta)$.
We may regard the function $\frac{d^2}{d\theta^2}f_\theta \Big{|}_{S_\theta}$ as a function, say $\tilde{f}$, on $V_\theta$ by pulling it back by $c^*\text{-Exp}_{y_\theta}$.
We apply Lemma~\ref{L:double double f_theta},
then under {\bf A3W} and the orthogonality \eqref{orthogonal to S_theta}, we see that the function $\tilde{f}$ is convex.

Now the $c^*$-convexity of $\Omega$ at each $y_\theta$ implies that
$V_\theta$ is convex, and so $\tilde{f}$ is a convex function on a convex domain.
Then we see by \eqref{gradient double f_theta at x_m}
that the function $\frac{d^2}{d\theta^2} f_\theta \Big{|}_{S_\theta}$ has the minimum value $0$ at $x_m$ so implying (1), and this minimum is unique under {\bf A3S} when $\tilde{f}$ is strictly convex so implying (2).
\end{proof}

\begin{remark}\label{R:weaker convex}
The proof of Proposition~\ref{P:double f positive} is the only place in the proof of Theorem~\ref{T:monotonicity}, where we use the $c^*$-convexity of $\Omega$ with respect to $\Lambda$. What we actually need here is that a dense subset of the set $(q_\theta + W_\theta ) \cap \text{dom}(c^*\text{-Exp}_{y_\theta})
$ can be reached by line segments from $q_\theta$ in $\text{dom}(c^*\text{-Exp}_{y_\theta})$. We can use this slight weaker condition, to deal with such examples as
$c(x,y) = -\frac{1}{2}\log |x-y|^2 $ in $\mathbf{R}^n \times \mathbf{R}^n \setminus \{ x=y\}$ (Example~\ref{E:log distance}); $c(x,y)=\pm \frac{1}{p}|x-y|^p, p \ne 0$ for $p = -2$ or $-\frac{1}{2} \le p < 1 $ ($-$ only) (Example~\ref{E:radial});
the cost $c(x,y)= -\frac{1}{2}\log |x-y|^2$ restricted on $S^{n-1} \subset \mathbf{R}^n$ which appears in the reflector antenna problem (Example~\ref{E:antenna}).
 \end{remark}

In the rest of the section, we prove Lemma~\ref{L:double double f_theta}.

\begin{proof}[Proof of Lemma~\ref{L:double double f_theta} {\bf (Sliding mountain and $c$-curvature)}]
We will use the same notation of Lemma~\ref{L:double double f_theta} and perform our calculations in the geodesic normal coordinates of $x_w$ and of $y_\theta$.
   Differentiating \eqref{parameter w} twice with respect to $w$ we see that
\begin{align}\label{zero x_w}
0= -D^2_{xx} D_y c(x_w , y_\theta) D_w x_w \ D_w x_w - D^2_{yx} c( x_w , y_\theta) D_{ww}^2 x_w  .
\end{align}

We compute by using \eqref{zero y_theta} and \eqref{zero x_w}
\begin{align}\nonumber
& D_{ww}^2 \frac{d^2}{d\theta^2} f_\theta (x_w) \\  \nonumber
 & = -D_{xx}^2 D_{yy}^2 c(x_w, y_\theta) D_w x_w  \ D_w x_w \ \dot{y}_\theta \cdot \dot{y}_\theta \\ \nonumber
& \ \ \ - D_x D_{yy}^2 c(x_w, y_\theta) D^2_{ww} x_w \ \dot{y_\theta} \cdot \dot{y_\theta} \\ \nonumber
& \ \ \
+ [-D^2_{xx} D_y c(x_w, y_\theta)  D_w x_w  \ D_w x_w- D^2_{yx}  c(x_w , y_\theta) D^2_{ww} x_w ] \cdot \ddot{y}_\theta  \\
\label{double double f_theta}
& = -D_{xx}^2 D_{yy}^2 c(x_w, y_\theta) D_w x_w  \ D_w x_w \ \dot{y}_\theta \cdot \dot{y}_\theta \\ \nonumber
& \ \ \ - D_x D_{yy}^2 c(x_w , y_\theta)
\Big{[} (D^2_{y, x} c(x_w , y_\theta) )^{-1} (-D^2_{xx}D_y c(x_w, y_\theta) D_w x_w \ D_w x_w \Big{]} \ \dot{y_\theta} \cdot \dot{y_\theta}, \end{align}

To compute $\mathfrak{S_c}(x_w , y_\theta) (D_w x_w , \dot{p}_\theta (x_w)  )$, we consider the line segment $p_\theta (x_w) + \psi \ \dot{p}_\theta (x_w)$, $ \psi \in \mathbf{R}$, in $T_{x_w} \Omega$, and its $c$-exponential image $y(\psi)$ as
\begin{align}\label{y_psi}
p_\theta (x_w) + \psi \ \dot{p}_\theta (x_w)
= -D_x c(x_w , y(\psi) ).
\end{align}
Note that $y(\psi) \big{|}_{\psi=0} = y_\theta$, and by differentiating \eqref{y_psi} with respect to $\psi$ that \begin{align*}
 \dot{p}_\theta (x_w)=
-D^2_{xy} c(x_w, y_\theta ) \frac{d}{d\psi} y(\psi) \Big{|}_{\psi=0},
 \end{align*}
therefore, we see
\begin{align*}
\dot{y}_\theta = \frac{d}{d\psi} y(\psi) \Big{|}_{\psi=0}
\end{align*} and by differentiating \eqref{y_psi} twice with respect to $\psi$ that
\begin{align}\label{zero y_psi}
 0= -D_{yy}^2 D_x c(x_w, y(\psi) )  \ \dot{y}_\theta  \ \dot{y}_\theta
 - D^2_{xy} c(x_w, y(\psi)) \frac{d^2}{d\psi^2} y(\psi) \Big{|}_{\psi =0}
 \end{align}

Now we compute using \eqref{zero y_psi}
\begin{align}\nonumber
 & \mathfrak{S_c}(x_w , y_\theta) (D_w x_w , \dot{p}_\theta (x_w)  ) \\ \nonumber
 & \ \ =  \frac{d^2}{d\psi^2}\Big{|}_{\psi =0}
 -D_{xx}^2 c(x_w, y(\psi)) D_w x_w \ D_w x_w  \\ \nonumber
& \ \ =-D^2_{yy} D^2_{xx} (x_w , y_\theta)  \ D_w x_w \ D_w x_w  \ \dot{y}_\theta \ \dot{y}_\theta\\ \nonumber
& \ \ \ \
- D_y D^2_{xx} c(x_w , y_\theta)   \
D_w x_w \ D_w x_w  \ \frac{d^2}{d\psi^2} y(\psi) \Big{|}_{\psi =0}\\
\label{S_c at x_w}
 & =  -D^2_{yy} D^2_{xx} (x_w , y_\theta) \ D_w x_w \ D_w x_w \ \dot{y}_\theta \ \dot{y}_\theta  \\ \nonumber & \ \ \ \
- D_y D^2_{xx} c(x_w , y_\theta)   \  D_w x_w \ D_w x_w
\Big{[}(D^2_{x y} c(x_w , y_\theta))^{-1}
(-D^2_{yy} D_x c(x_w , y_\theta) \ \dot{y}_\theta \ \dot{y}_\theta ) \Big{]}
\end{align}

Therefore, comparing \eqref{double double f_theta} and \eqref{S_c at x_w} by interchanging the order of the derivatives and by Lemma~\ref{L:symmetry}, we see
that
\begin{align*}
 & D_{ww}^2 \frac{d^2}{d\theta^2} f_\theta (x_w)
 =\mathfrak{S_c}(x_w , y_\theta) (D_w x_w , \dot{p}_\theta (x_w)  ) .
\end{align*}
This finishes the proof of Lemma~\ref{L:double double f_theta}.
\end{proof}

\section{Examples}\label{S:examples}
In the following, we demonstrate our result in some examples of cost functions among the ones
given by Ma, Trudinger, and Wang \cite{MTW} \cite{TW}, and by Loeper \cite{L}.  A common
feature of these examples is that they only require the $c$-convexity of the target domain
$\Lambda$ with respect to the source domain $\Omega$, not the other direction,  for the
equivalence $\partial^c \phi =\partial \phi$ of $c$-subdifferential and ordinary
subdifferential of $c$-convex functions $\phi$ on $\Omega$ indexed by $\Lambda$.

 \begin{example}\label{E:convex f g}
$c(x,y) = |x-y|^2 + |f(x)-g(y)|^2$
satisfies {\bf A0, A1, A2, A3W} (resp. {\bf A3S}) for $\Omega=\Lambda=\mathbf{R}^n$, if $f, g :\mathbf{R}^n \to \mathbf{R}$ convex (resp. strictly convex)
and $|\nabla f|, |\nabla g| <1$.
Such $\Omega$ and $\Lambda$ are $c$, $c^*$-convex with respect to the other, since the domains of $c$, $c^*$ exponential maps are the whole space $\mathbf{R}^n$ (One can verify this for example, by showing that the map $-\nabla_x c (x,\cdot)$ is a proper map.) So by Corollary ~\ref{C:subdomains}, we see that for any domains $\Omega_1$, $\Lambda_1$ $\subset \mathbf{R}^n$, $\Lambda_1$ being $c$-convex with respect to $\Omega_1$, and for any $c$-convex function $\phi$ on $\Omega_1$ indexed by $\Lambda_1$,
we have $\partial^c \phi = \partial \phi$ without the $c^*$-convexity of $\Omega_1$.
This result is in fact easy to show by directly computing
$\frac{d^2}{d\theta^2}f_\theta$ and checking its positivity.
\end{example}

\begin{example}\label{E:sqrt}
We can conclude the same as in Example~\ref{E:convex f g} for the cost function $c(x,y)=\sqrt{1+ |x-y|^2}$ on $\mathbf{R}^n \times \mathbf{R}^n$.
\end{example}

In the following examples, we will use the notation of Section~\ref{S:Proof T:c-subdifferential}. These examples concern cost functions that have singularities.

\begin{example}\label{E:log distance}
$c(x,y) = -\frac{1}{2} \log |x-y|^2$ on
$\mathbf{R}^n \times \mathbf{R}^n \setminus \{ x=y\}$ satisfies {\bf A3S}. In this example, one can show that the curve
$\{ y_\theta\}_{\theta \in \mathbf{R}}$ is a circle (possibly with infinite radius) passing through $x_m$, tangent to $p_1 -p_0$ at $x_m$, and that
$S_\theta$'s are $(n-1)$-dimensional spheres (with different radius and centers, possibly with infinite radius) passing through $y_\theta$ and $x_m$, with normal vector $p_1 - p_0$ at $x_m$. This geometric configuration shows directly the monotonicity
of $S_\theta^+$, as in Theorem~\ref{T:monotonicity}, and thus either for $\Omega=\Lambda=\mathbf{R}^n$ or for any $\Omega \subset \mathbf{R}^n$, and  $\Lambda \subset \mathbf{R}^n$, $c$-convex with respect to $\Omega$, we have
$\partial^c \phi = \partial \phi$ for all $c$-convex functions $\phi$ on $\Omega$ indexed by $\Lambda$.  The same result can be obtained applying our general theory, especially Theorem~\ref{T:monotonicity} and  Remark~\ref{R:weaker convex}.  Note that in this case
the domains $\text{dom}(c\text{-Exp}_x)$, $\text{dom}(c^*\text{-Exp}_y)$ (for $\Omega=\Lambda=\mathbf{R}^n$) are $\mathbf{R}^n \setminus \{0\}$ and so
\begin{align*}
(q_\theta + W_\theta ) \cap \text{dom}(c^*\text{-Exp}_{y_\theta})=(q_\theta + W_\theta ) \setminus \{0\},
\end{align*}
(see \eqref{S_theta}\eqref{W theta} for the definition of $W_\theta$).
The set $(q_\theta + W_\theta ) \setminus \{0\}$ is convex except when $\{0\} \subset q_\theta + W_\theta$.
In the latter case, still a dense subset of it can be reached by line segments from $q_\theta$ within it,
  if the dimension $n \ge 3$, as we needed in Remark~\ref{R:weaker convex}.
The case $n=1$ is trivial, and the case $n=2$ can be handled by an approximation argument.

One may be concerned with that when $p_0, 0, p_1$ are collinear in the tangent space at $x_m$, the corresponding $c$-segment is no longer
differentiable when $p_\theta$ passes through the origin $0$ since the point $y_\theta$
in the $c$-segment goes to infinity. This problem can be bypassed by approximating the segment
$[p_0, p_1]$ by a line segment which does not pass through the origin. Then our desired
{\bf DASM} survives under the approximation, and so it holds for $\Omega=\Lambda=\mathbf{R}^n$.
\end{example}

\begin{example}\label{E:radial}

$c(x,y)=\pm \frac{1}{p}|x-y|^p, p \ne 0$ satisfies {\bf A3W} for $p = \pm 2$ or $p=-\frac{1}{2}$ ($-$ only) and {\bf A3S} for $-\frac{1}{2} < p < 1 $ ($-$ only).
The case $p=2$ is well known and the same conclusion as in Example~\ref{E:convex f g} holds. For other $p$, the cost is not differentiable for $x=y$, but this case is similar to Example~\ref{E:log distance} and by the same approximation method using Remark~\ref{R:weaker convex} we get the same conclusion.
\end{example}

\begin{example}\label{E:sphere}{\bf (The round sphere)} We consider the cost
$c(x,y) = \frac{1}{2}d^2 (x,y)$ on $S^n$, with $d$ the Riemannian distance function of the standard round metric of $S^n$, i.e. with the diameter $\pi$. Loeper showed \cite{L} that
$c$ satisfies {\bf A3S} by directly computing the cost-sectional curvature $\mathfrak{C}_s$.
In this case, $\text{dom}(c^*\text{-Exp}_{y_\theta})= B_{\pi} (0) \subset T_{y_\theta} S^n$ (for $\Omega =S^n$) and so
$(q_\theta + W_\theta ) \cap \text{dom}(c^*\text{-Exp}_{y_\theta})$ is convex (see \eqref{S_theta}\eqref{W theta} for the definition of $W_\theta$) and we can apply Proposition~\ref{P:double f positive}
 for those points where the function $c(\cdot, y_\theta)$ is differentiable, to get {\bf Strict DASM} on the set $S^n \setminus \cup_{0 \le \theta \le 1} \hat{y}_{\theta}$, where we denote by $\hat{y}_\theta$ the antipodal point of $y_\theta$.
Now, we use the continuity of the distance function to extend the {\bf Strict DASM} to the whole $S^n$. From this result together with Remark~\ref{R:double-sliding}, we can conclude that  if either $\Omega=\Lambda=S^n$ or $\Lambda \subset S^n$ is $c$-convex with respect to $\Omega \subset S^n$ then  $\partial^c \phi = \partial \phi$ for all $c$-convex functions $\phi$ on $\Omega$ indexed by $\Lambda$.


 \end{example}

\begin{example}\label{E:antenna}
{\bf (Reflector antenna problem)}
The cost $c(x,y)= -\frac{1}{2}\log |x-y|^2$ restricted on the unit sphere $S^{n-1} \subset \mathbf{R}^n$ satisfies {\bf A3S} on $S^{n-1} \times S^{n-1} \setminus \{ x=y\}$. In this case, for $\Omega= S^{n-1}$,
\begin{align*}
\text{dom}(c^*\text{-Exp}_{y})=  T_{y} S^{n-1} \cong \mathbf{R}^{n-1}
\end{align*}
and $c^*\text{-Exp}_{y}(T_{y} S^{n-1})= S^{n-1} \setminus \{y\}$, $c^*\text{-Exp}_{y}(0)= \hat{y}$, where we denote by $\hat{y}$ the antipodal point of $y$ in $S^{n-1}$. So by our general theory (Theorem~\ref{T:monotonicity} and  Remark~\ref{R:weaker convex})
as applied in Example~\ref{E:log distance},
 we can make the conclusion that if either $\Omega=\Lambda=S^n$ or $\Lambda \subset S^n$ is $c$-convex with respect to $\Omega \subset S^n$, then  $\partial^c \phi = \partial \phi$ for all $c$-convex functions $\phi$ on $\Omega$ indexed by $\Lambda$.

\end{example}

 \bibliographystyle{plain}

\end{document}